\documentclass[11pt,twoside]{article}

\usepackage{a4wide}
\usepackage{amssymb}
\usepackage[matrix,arrow,curve]{xy}
\usepackage{epsfig}

\usepackage{amsmath}
\usepackage{color}
\usepackage{hyperref}
 \definecolor{darkblue}{rgb}{0,0,.7}
 \hypersetup{pdftex=true, colorlinks=true, breaklinks=true, linkcolor=darkblue, citecolor=darkblue, menucolor=darkblue, pagecolor=blue, urlcolor=darkblue}

\numberwithin{equation}{section}
\numberwithin{figure}{section}

\usepackage{theorem}
\newtheorem{thm}{Theorem}[section]
\newtheorem{lem}[thm]{Lemma}

\newtheorem{cor}[thm]{Corollary}

\newtheorem{remark}[thm]{Remark}

\newtheorem{defi}[thm]{Definition}

\newcommand{\comment}[1]{}

\newcommand{\Hom}{{\text{Hom}}}

\newcommand{\Z}{\mathbb Z}

\newcommand{\Pro}{\mathrm{P}^{\infty}}
\newcommand{\Po}{{\mathrm{P}}}
\newcommand{\HB}{\mathbb H}

\newcommand{\R}{\mathbb R}
\newcommand{\Q}{\mathbb Q}
\newcommand{\C}{\mathbb C}
\newcommand{\qed}{\quad \textbf{q.e.d.}}

\DeclareMathOperator{\Tor}{Tor}
\DeclareMathOperator{\Ext}{Ext}

\DeclareMathOperator{\Kern}{Ker}
\DeclareMathOperator{\Coker}{Coker}

\xyoption{frame}
\xyoption{tips}
\xyoption{color}
\entrymodifiers={+!!<0pt,\fontdimen22\textfont2>}

\title{Conjugations on 6-manifolds with free integral cohomology}
\author{Martin Olbermann}
\date{\today}
\begin{document}
\maketitle

\begin{abstract}
\noindent In this article, we show the existence of conjugations on many simply-connected spin 6-manifolds with free integral cohomology.
In a certain class the only condition on $X^6$ to admit a conjugation with fixed point set $M^3$
is the obvious one: the existence of a degree-halving ring isomorphism between the $\Z_2$-cohomologies of $X$ and $M$.
\end{abstract}

%\tableofcontents

\section{Introduction}

A conjugation space is a space together with an involution satisfying 
a certain cohomological pattern. The definition of a conjugation space \cite{HHP} 
has its origin in the observation that there are many examples of involutions with this cohomological
pattern, e.g. the complex conjugation on the projective space $\C \Po^n$ and on
complex Grassmannians, natural involutions on smooth toric manifolds \cite{DJ} and on 
polygon spaces \cite{HK}. There is always a degree-halving isomorphism $\kappa$ between 
the $\Z_2$-cohomologies of the space and its fixed point set.
Hausmann, Holm and Puppe \cite{HHP} also discovered that in all these examples there is an even richer
structure in the equivariant cohomology of the space itself and the space of fixed points 
of the involution. This extra structure is also part of the definition of a conjugation space.
Their main examples of conjugation spaces are conjugation complexes:
every cell complex with the property that each cell is a unit disk in $\C^n$ with complex conjugation, 
and with equivariant attaching maps, is a conjugation space. 
They also prove that coadjoint orbits of semi-simple Lie groups with the Chevalley involution are conjugation spaces.
One can form connected sums of manifolds with conjugations to obtain new conjugation spaces.
Moreover, there are several bundle constructions yielding conjugations on the total space, 
provided conjugations on base and fiber.

\bigskip

As a general assumption, we are interested in smooth conjugations on connected, closed, smooth manifolds.
These have to be of even dimension.
A point with the trivial involution is a conjugation space.
In dimension 2, only the sphere admits conjugations.
Dimension 4 is special since the quotient space of the involution can be given a smooth structure \cite{HH}. 
Hambleton and Hausmann prove that equivariant diffeomorphism classes of oriented connected conjugation 
4-manifolds correspond to diffeomorphism classes of pairs $(X,\Sigma)$, where $X$ 
is an oriented 4-dimensional $\Z_2$-homology sphere and $\Sigma$ is a closed connected subsurface of $X$.
Hence the study of conjugations is ``reduced" to a difficult non-equivariant problem: a generalization of
four-dimensional knot theory. 
Another result by Hambleton and Hausmann is that a simply-connected spin 4-manifold 
admitting a conjugation is homeomorphic a connected sum of $S^2\times S^2$'s. Thus there are for example no smooth conjugations
on the Kummer surface, although its $\Z_2$-cohomology ring is the same as the one of a connected sum of $S^2\times S^2$'s.
 
The problem of realizing 3-manifolds as fixed point sets of conjugations was considered in \cite{OlbThesis}.
The present paper focuses on the identification of the resulting conjugation space. It shows
that for a certain class of 6-manifolds, the obvious condition on $X^6$ to admit a conjugation with fixed point
set a given $M^3$ (namely the existence of a degree-halving isomorphism between the $\Z_2$-cohomologies)
is the only condition.

\bigskip

We focus on a class of 6-manifolds which has been considered in connection with the question of the existence
of simply-connected asymmetric manifolds. 
In the smooth category, a manifold $X$ is called asymmetric if every (smooth) action of a finite (or equivalently compact Lie) group on $X$ is trivial.
One of the problems in finding simply-connected closed asymmetric manifolds is that every explicit construction of a manifold
seems to produce a non-trivial symmetry group. Thus we should be interested in more implicit descriptions of the manifolds under
consideration. For example we could use abstract classification results as the following, proved by Wall. 

\begin{thm}[\cite{Wal}] \label{WallThm}
Diffeomorphism classes of six-dimensional manifolds which are simply-connected, spin, have free integer cohomology,
and zero odd cohomology, correspond bijectively to isomorphism classes of the following invariants:
\begin{enumerate}
\item
the rank $m$ of the second cohomology  with integer coefficients $H^2$,
\item
a trilinear form $\mu$ on $H^2$ corresponding to the evaluation of triple cup products on the fundamental class,
\item
and a linear form $P_1$ on $H^2$ corresponding to the Poincar\'e dual of the first Pontrjagin class, 
\end{enumerate}
subject to the conditions:
\begin{eqnarray*}
\mu(x,x,y)&\equiv& \mu(x,y,y) \ \text{(mod 2)}\ \text{ for all } x,y\in H^2, \\
P_1(x)&\equiv& 4\mu(x,x,x) \ \text{(mod 24)}\ \text{ for all }x \in H^2.
\end{eqnarray*}
\end{thm}

V.~Puppe considered the class of such manifolds and proved certain asymmetry statements for 
manifolds in this class \cite{P6}. In particular Puppe proved the following theorem which gives a 
relation between conjugations and asymmetric manifolds.

\begin{thm}[\cite{P6}] \label{Pthm}   
Some simply-connected spin 6-manifolds do not admit a
non-trivial \textbf{orientation-preserving} $\Z_p$-action.
Moreover, they are either asymmetric, or conjugation spaces.
\end{thm}

%(By a recent result of Kreck \cite{KrA}, conjugations are not possible on infinitely many of Puppe's examples. Hence they are asymmetric manifolds.)

\bigskip

In our construction \cite{OlbThesis} of 6-dimensional manifolds with conjugations we (necessarily) obtain manifolds 
whose $\Z_2$-cohomology is concentrated in even degrees. 
This implies that for integral cohomology, no torsion of even order is possible. 
Odd torsion might occur, however: in the last step of the construction (surgery in the middle dimension), 
we only control the $\Z_2$-cohomology of the result, since we find a hyperbolic summand only in $\Z_2$-cohomology.
In the first part of this article, a more detailed analysis of the intersection form in the middle dimension leads to a better
choice for the surgery in the middle dimension. (In particular we use information given by universal coefficient spectral sequences.)
From this we deduce our first result:

\begin{thm}\label{mthm}
Every closed orientable 3-dimensional manifold such that $H_1(M;\Z)$ contains no 2-primary torsion can be realized as fixed point space of a smooth conjugation 
on a closed simply-connected spin 6-manifold with free integral and zero odd cohomology.
\end{thm}
(Under the conjugation space isomorphism $\kappa$ of $\Z_2$-cohomology algebras 
dividing the degree by two, the second Wu class of the 6-manifold will be mapped to the first Wu class of its fixed point set, 
since we have $\kappa(Sq^{2k}(x))=Sq^k(\kappa(x))$, see \cite{FP}. This implies that only orientable 3-manifolds may occur as fixed point sets of conjugations on such 6-manifolds.)

\bigskip

In the second part, we analyze which 6-manifolds are produced by our construction. 
The main result of the article is
\begin{thm}\label{nthm}
Let $X$ be a closed simply connected spin 6-manifold with free and only even integral cohomology.
Let $M$ be a closed oriented 3-manifold such that $H_1(M;\Z)$ contains no 2-primary torsion, and such that there exists a ring isomorphism
$H^*(X;\Z_2)\to H^*(M;\Z_2)$ dividing all degrees by 2. Then there exists a smooth conjugation on
$X$ with fixed point set $M$.
\end{thm}
From this we deduce:
\begin{cor}\label{mcor}
Every closed simply connected spin 6-manifold with free and only even integral cohomology, such that
$x^2$ is divisible by 2 for all $x\in H^2(M;\Z)$ admits a conjugation.
In particular the examples from \cite{P6} are not asymmetric, but do admit conjugations.
\end{cor}

\begin{remark}
This corrects the wrong statement in \cite{KrA} that Puppe's examples are asymmetric. In order to modify his argument such that
it gives asymmetric manifolds, Kreck moved to the non-smooth category \cite{KrC}. 
\end{remark}
Another immediate corollary is a partial converse to P.A. Smith's theorem about group actions on spheres.
\begin{cor}
Every $\Z_2$-homology 3-sphere is the fixed point set of an involution (a conjugation) on $S^6$.
\end{cor} 

In a forthcoming paper \cite{sphinv}, we show that our surgery approach to conjugations, 
and the computation of the transfer map, can be used not only for existence questions,
but are also suited to give classification results.
In particular we classify all involutions on $S^6$ with three-dimensional fixed point set.

\bigskip

The obvious open question is whether in the main theorem \ref{nthm} the condition that
$x^2$ is divisible by 2 for all $x\in H^2(M;\Z)$ is necessary. An example where the condition is violated is that
every homotopy projective $\C \Po^3$ admits a conjugation with fixed point set
$\R \Po^3$. 
(Dovermann, Masuda and Schultz show the existence of such involutions \cite{DMS}, 
and all these involutions must be conjugations.) 

\bigskip

\noindent {\bf{Acknowledgements.}} The author was partially supported by a research fellowship of the DFG. 
We would like to thank Diarmuid Crowley, Ian Hambleton, Jean-Claude Hausmann, Matthias Kreck and Arturo Prat-Waldron 
for many helpful discussions and remarks.

\section{Construction of conjugation 6-manifolds with free integral cohomology}
\subsection{Conjugations on manifolds}
A conjugation on a topological space $X$ is an involution $\tau:X\to X$, which we consider as an action of the group $\Z_2\cong C= \{id, \tau\}$
on $X$, and which satisfies the following cohomological pattern:
We denote the Borel equivariant cohomology of $X$ by $H^*_C(X;\Z_2)$. It is a module over $H^*_C(pt;\Z_2)=\Z_2[u]$. 
The restriction maps in equivariant cohomology are denoted by $\rho: H^*_C(X;\Z_2)\to H^*(X;\Z_2)$
and $r: H^*_C(X;\Z_2)\to H^*_C(X^\tau;\Z_2)\cong H^*(X^\tau;\Z_2)[u]$.

\begin{defi}
$X$ is a conjugation space if 
\begin{itemize}
\item
$H^{odd}(X;\Z_2)=0$, 
\item there exists a (ring) isomorphism $\kappa: H^{2*}(X;\Z_2)\to H^*(X^\tau;\Z_2)$
\item and a (multiplicative) section $\sigma: H^*(X;\Z_2)\to H^*_C(X;\Z_2)$ of $\rho$ 
\item such that the so-called conjugation equation holds:
$$r\sigma(x)=\kappa(x)u^k + \text{ terms of lower degree in }u.$$
\end{itemize}
\end{defi}
One does not need to require that $\kappa$ and $\sigma$ be ring homomorphisms, it is a consequence of the definition. Moreover,
the ``structure maps" $\kappa$ and $\sigma$ are unique, and natural with respect to equivariant maps between conjugation spaces.

A conjugation manifold is a conjugation space consisting of a smooth manifold $X$ with a smooth involution $\tau$.
As a consequence, $X$ must be even-dimensional, say of dimension $2n$, and $M$ is of dimension $n$. 

In \cite{OlbThesis} we proved that it is possible to give a definition of conjugation spaces without
the non-geometric maps $\kappa$ and $\sigma$, which is moreover well-adapted to the case of conjugation manifolds,
where the fixed point set has an equivariant tubular neighbourhood. An application of this characterization is 
theorem \ref{Wcond} below.

\subsection{Review and outline}

In this section we review the construction from \cite{OlbThesis}, and we explain how it is modified in the next section.

\begin{thm}[\cite{OlbThesis}] \label{Wcond} 
Let $X$ be a closed 6-manifold with a smooth involution $\tau$
such that the fixed point set $M$ is a 3-manifold with trivial normal
bundle. Using the equivariant tubular neighbourhood theorem, we write
$X=(M\times D^3) \cup V$, where the involution restricts to a free involution $\tau$ on $V$ such that $W:=V/ \tau$ is a
$6$-manifold with boundary $\partial W=M\times \R \Po^2$.
Then $X$ is a conjugation space if and only if restriction to the boundary $$H^*(W;\Z_2)\to H^*(M\times \R \Po^2;\Z_2)$$ induces an isomorphism:
$$H^*(W;\Z_2)\to H^*(M\times \R \Po^2;\Z_2) / \bigoplus_{i>j} H^i(M;\Z_2) u^j.$$
Translated to homology this is equivalent to the condition that inclusion of the boundary
$$H_*(M\times \R \Po^2;\Z_2)\to H_*(W;\Z_2)$$ induces an
isomorphism: $$\bigoplus_{i\le j} H_i(M;\Z_2)\otimes H_j(\R \Po^2;\Z_2) \to H_*(W;\Z_2).$$
\end{thm}
Thus, in order to construct a conjugation on some 6-manifold with fixed point set $M$, it suffices to
find the ``right" nullbordism $W$ of $M\times \R \Po^2$. (The trivial normal bundle condition is always
satisfied if $M$ is oriented and $X$ is simply-connected.)

We want the manifolds $X$ to be among those classified by Wall, 
which means simply-connected spin and with $H_3(X)=0$. This gives conditions on $W$, and on its normal 2-type: 

\begin{defi}
The normal 2-type of a compact manifold $N$ is a fibration $B_2(N)\to BO$ which is obtained
as a Postnikov factorization of the stable normal bundle map $N\to BO$:
There is a 3-connected map $N\to B_2(N)$, the fibration 
$B_2(N)\to BO$ is 3-coconnected (i.e. $\pi_i(BO,B_2(N))=0$ for $i>3$), and
the composition is the stable normal bundle map. This determines $B_2(N)\to BO$ up to fiber homotopy equivalence.
\end{defi}

Given a closed oriented 3-manifold $M$, we set $m=rk(H_1(M;\Z_2))+1$. 
We construct a fibration $B\to BO$ which is 3-coconnected, and such that $B$ is connected, $\pi_1(B)=\Z_2$,
$\pi_2(B)\cong \Z^m$, on which $\pi_1(B)$ acts by multiplication with $-1$, and $\pi_3(B)=0$.
More precisely, we use the fiber bundle $S^2\to \C \Pro \to \HB \Pro$ induced from the identification $\C^{\infty}\cong \HB^{\infty}$.
The antipodal map on each fiber induces a free involution $\tau$ on $\C \Pro$. 
Now $B=BSpin\times Q_m$, where $$Q_m=(\C \Pro \times \dots \times \C \Pro \times S^{\infty})/(\tau , \dots , \tau, -1),$$
where we have $m$ factors $\C \Pro$. We have the projection on the last factor $p:Q_m\to \R\Pro=BO(1)$.
Since the involution $\tau$ is free, we can identify up to homotopy $Q_m\cong (\C \Pro)^m/\tau^m$.
The map $B\to BO$ is the composition $B=BSpin\times Q_m \stackrel{B\pi\times
p}{\to} BO\times BO(1) \stackrel{\oplus}{\to} BO$.

\begin{lem}[\cite{OlbThesis}]
If $X$ is simply-connected spin, with free and only even cohomology, $\tau$ is a conjugation on $X$ with fixed point set $M$ 
which acts by $-1$ on $H^2(X)$, and $W$ is defined as before, then the normal 2-type of $W$ is $B\to BO$.
\end{lem}

Now start with a closed oriented 3-manifold $M$.
By a bordism calculation, we show that there exists a manifold $W$ with boundary $M\times \R \Po^2$ and normal 2-type $B\to BO$:
\begin{thm}[\cite{OlbThesis}]
The bordism group of 5-manifolds with normal $B$-structures is trivial: $\Omega_5^B=0$.
(A normal $B$-structure is a lift of the normal bundle map to $B$. For a more precise definition, see
\cite{Sto}.)
Using surgery below the middle dimension, we may assume that $B\to BO$ is the normal 2-type for $W$.
\end{thm}

Thus the map $W\to B$ is 3-connected, and so in particular $\pi_1(W)=\Z_2$. 
We would like to obtain a manifold $W_0$ with the same boundary, and such that we get an isomorphism
$$\bigoplus_{i\le j} H_i(M;\Z_2)\otimes H_j(\R \Po^2;\Z_2) \cong H_*(W_0;\Z_2).$$ 
For $W$ this map is an isomorphism except in dimension 3, where the map is injective, but has a possibly non-trivial cokernel. 
It is this cokernel we want to kill using surgery. The cokernel is isomorphic to $H_3(W;\Z_2)/rad$, where $rad$ is the
radical of the intersection form on $H_3(W;\Z_2)$. We find that the intersection form on $H_3(W;\Z_2)/rad$ is hyperbolic, 
and we find disjointly embedded $3$-spheres in $W$ with trivial normal bundle mapping to generators for a Lagrangian.
We use the following lemma.
\begin{lem}[\cite{OlbThesis}]
The Hurewicz map $\pi_3(W)\to H_3(W;\Lambda)$ and the map $H_3(W;\Lambda)\to H_3(W;\Z_2)/rad$ are both surjective.
\end{lem}
By surgery on these $3$-spheres we kill the cokernel in the middle dimension and obtain a manifold $W_0$ with the desired properties.
The double cover $V_0$ of $W_0$ has boundary $M\times S^2$. We glue $V_0$ along its boundary to $M\times D^3$
and obtain a 6-manifold $X=M\times D^3 \cup V_0$. The involution on $X$ which is $(id,-id)$ on $M\times D^3$ and which equals the unique 
non-trivial deck transformation on $V_0$ is a conjugation with fixed point set $M$.

\begin{thm}[\cite{OlbThesis}]
Every closed orientable 3-manifold can be realized as fixed point set of a smooth conjugation on a closed simply-connected
spin 6-manifold.
\end{thm}

The drawback of the surgery procedure is that we have not controlled the $\Lambda=\Z [\Z_2]$-valued quadratic form on $H_3(W;\Lambda)$.
(We only extracted partial information which made sure that we found disjointly embedded $3$-spheres in $W$ 
with trivial normal bundle we could do surgery on.) So we cannot say precisely which 6-manifolds with conjugations we obtain.

\bigskip

In the following we will show that $H_3(W;\Lambda)/rad$ is a free $\Lambda$-module and that the induced quadratic form 
is non-degenerate and thus stably hyperbolic.
(Morally, the quadratic form is given by the map $H_3(W;\Lambda)\to H_3(W,\partial W;\Lambda)$, and the absence of (odd) $\Z$-torsion in 
$H_2(\partial W;\Lambda)$ implies that it is possible to do the necessary surgeries without creating odd torsion.)
The intersection form on $H_3(W;\Lambda)$ maps to the $\Z_2$-valued intersection form on $H_3(W;\Z_2)$ 
using the map $\epsilon:\Lambda\to \Z_2$ defined by $a+bT\mapsto a+b$.

We conclude that (after possibly stabilizing $W$ by connected sum with copies of $S^3\times S^3$) we can do surgery on generators of a Lagrangian for $H_3(W;\Lambda)/rad$. 
This kills $H_3(W;\Z_2)/rad$ as before, but we also know that the effect on $\Lambda$-homology is just to kill $H_3(W;\Lambda)/rad$. 
In particular we see that the normal 2-type of $W_0$ is equal to the normal $2$-type of $W$. By the Mayer-Vietoris sequence for $X=M\times D^3 \cup V_0$,
we see that $H_2(X)$ is a free $\Z$-module on which the conjugation acts by multiplication with $-1$ (see \cite{OlbThesis}).
Since $H_1(X)$ and $H_2(X)$ are free over $\Z$, Poincar\'e duality implies that all homology of $X$ is free over $\Z$.
But since we have a degree-halving $\Z_2$-homology isomorphism from $X$ to $M$, we see that the free homology of $X$ 
is concentrated in even degrees.

\subsection{Proof of theorem \ref{mthm}}\label{ConCon}

Recall from \cite{OlbThesis} that the normal 2-type for $W$ is $B=BSpin\times Q_m \to BO$, where 
$Q_m$ is a fibration over $\R \Po^\infty$ with fiber $K(\Z^m,2)$, and $\pi_1(B)$ acts on $\pi_2(B)$ by multiplication with $-1$.
Since $\pi_3(B)=0$ we have a $(-1)$-quadratic form on $\pi_3(W)=\Kern(\pi_3(W)\to \pi_3(B))$.  
Let $\Lambda=\Z [\Z_2 ]\cong \Z[T]/\langle T^2-1\rangle $, with the involution $a+bT\mapsto a-bT$. Then the quadratic form consists of a 
$(-1)$-hermitian map $\lambda: \pi_3(W)\times \pi_3(W)\to \Lambda$ and a quadratic refinement $\tilde{\mu}:\pi_3(W)\to \Lambda / \Z\cdot 1$.
(This is a slight modification of the usual quadratic form \cite{SCM} counting (self-)intersections. 
Since we consider the self-intersection on elements of the homotopy group and not on regular homotopy classes of immersions, 
the values of $\tilde{\mu}$ are defined only modulo $\Z\cdot 1$. See \cite{Kre} for details.)

\bigskip

The map $\lambda$ fits in a commutative diagram with the Hurewicz map
$\pi_3(W)\to H_3(W;\Lambda)$ and the $\Lambda$-valued intersection
pairing on $H_3(W;\Lambda)$. 
It also maps to  $H_3(W;\Z_2)$ and the $\Z_2$-valued intersection pairing, 
using the map $\epsilon:\Lambda\to \Z_2$ defined by $a+bT\mapsto a+b$.

We have seen that the map $\pi_3(W)\to H_3(W;\Lambda)$ is surjective,
and that the map $H_3(W;\Lambda)\to H_3(W;\Z_2)$ induces a 
surjection onto $H_3(W;\Z_2)/rad$, and that the latter carries
a hyperbolic bilinear form.

\bigskip

It turns out that a detailed analysis of the long exact sequences of the pair $(W,\partial W)$ with $\Lambda$- and $\Z_2$-coefficients leads to a good understanding of the quadratic form. 
So let us analyze the commutative diagram 
$$ 
\xymatrix@C-17pt{
H_4(W,\partial W;\Lambda)
\ar[r]
\ar[d]
& H_3(\partial W;\Lambda)
\ar[d]
\ar[r]
& H_3(W;\Lambda)
\ar[r]
\ar[d]
& H_3(W,\partial W;\Lambda)
\ar[r]
\ar[d]
& H_2(\partial W;\Lambda)
\ar[d]
\ar[r]
& H_2(W;\Lambda)
\ar[d]
\\
H_4(W,\partial W;\Z_2)
\ar[r]
& H_3(\partial W;\Z_2)
\ar[r]
& H_3(W;\Z_2)
\ar[r]
& H_3(W,\partial W;\Z_2)
\ar[r]
& H_2(\partial W;\Z_2)
\ar[r]
& H_2(W;\Z_2)
}
$$

In the following we want to use (``universal") Poincar\'e duality with group ring coefficients.
For this we need oriented covers (see \cite{Ran}, Def.~4.56) of our manifolds.
The homology groups with group ring coefficients are just the integral homology groups of the oriented covering space, 
with left action of the fundamental group given by deck transformations.
The definition of cohomology groups with group ring coefficients is more difficult. For finite fundamental groups,
they are defined in the following way: The cohomology groups of the covering space 
are naturally right modules over the group ring, if one uses the action of the fundamental group given by deck transformations.
In order to make them left modules, we need an involution on the group ring. 
We use the ``involution twisted by the first Stiefel-Whitney class".
In our case this is $a+bT\mapsto a-bT$. The result is that the $\Lambda$-cohomology groups of our manifolds
are the integral cohomology groups of their double covers, but $T$ acts by $(-1)$ times the deck transformation.

With this definition we have universal Poincar\'e duality, i.e. isomorphisms of $\Lambda$-modules
\begin{gather*}
H_k(W;\Lambda)\cong H^{6-k}(W,\partial W;\Lambda),\\ 
H^k(W;\Lambda)\cong H_{6-k}(W,\partial W;\Lambda), \\
H_k(\partial W;\Lambda) \cong H^{5-k}(\partial W;\Lambda).
\end{gather*}
More precisely, universal Poincar\'e duality is described in \cite{Ran}.

Using Poincar\'e duality,
$H_4(W,\partial W; \Lambda)\to H_3(\partial W; \Lambda)$ 
identifies with the map $H^2(W;\Lambda)\to H^2(\partial W;\Lambda)$. 
We have $H^2(W;\Lambda)\cong \Z^m_+$ and 
$H^2(\partial W;\Lambda)\cong H^2(M;\Z)_-\oplus H^2(S^2)_+$
(we indicate the $T$-actions by subscripts $\pm$).

So we have to look at the double cover of the map $M\times \R \Po^2 \to Q_m$.
Going through our construction of this map \cite{OlbThesis}, we see that
the map $M\times S^2 \to (\C \Po^\infty)^m$ maps the 2-skeleton of $M$ to a point,
and $S^2=\C \Po^1$ diagonally into $(\C \Po^\infty)^m$.
Hence the map $$\Z^m_+\to H^2(M;\Z)_-\oplus H^2(S^2)_+$$ maps every $\Z$-summand 
isomorphically to $H^2(S^2)_+$, and has cokernel 
$H^2(M;\Z)\cong H_1(M)_+\otimes H_2(S^2)_-$.
And $H_2(\partial W;\Lambda)\to H_2(W;\Lambda)$ is the map 
$H_2(M)_+\oplus H_2(S^2)_-\to \Z^m_-$ sending $H_2(M)_+$ to 0 and the generator 
of $H_2(S^2)_-$ to $(1,1,\dots, 1)$.
We get an exact sequence:
$$0\to H_1(M)_+\otimes H_2(S^2)_- \to H_3(W;\Lambda) \to H_3(W,\partial W; \Lambda)\to H_2(M)_+\to 0.$$
Note that $H_2(M)_+$ is free over $\Z$. 

\bigskip

From the condition that ensures that we obtain a conjugation space (which is fulfilled for $W$ except in the middle dimension) we know that the map on $H_2$ with $\Z_2$-coefficients restricts to an isomorphism
$$H_0(M;\Z_2)\otimes H_2(\R \Po^2;\Z_2)\oplus 
H_1(M;\Z_2)\otimes H_1(\R \Po^2;\Z_2) \to H_2(W;\Z_2).$$
Thus the map $H_2(\partial W;\Z_2) \to H_2(W;\Z_2)$ has kernel isomorphic to $H_2(M;\Z_2)\cong \Z_2^{m-1}$. 
By Poincar\'e duality, the map $H_4(W,\partial W;\Z_2) \to H_3(\partial W;\Z_2)$ has cokernel 
$\Z_2^{m-1}$.
So we get the commutative diagram with exact rows:
$$ 
\xymatrix{
0
\ar[r]
& H_1(M)_+\otimes H_2(S^2)_- 
\ar[r]
& H_3(W;\Lambda)
\ar[r]
\ar[d]
& H_3(W,\partial W;\Lambda)
\ar[r]
\ar[d]
& H_2(M)_+
\ar[r]
& 0
\\
0
\ar[r]
& \Z_2^{m-1}
\ar[r]
& H_3(W;\Z_2)
\ar[r]
& H_3(W,\partial W;\Z_2)
\ar[r]
& \Z_2^{m-1}
\ar[r]
& 0
}
$$

\bigskip
By Poincar\'e duality we see that $H_3(W,\partial W; \Lambda)\cong H^3(W;\Lambda)$ 
is $\Z$-free since its $\Z$-torsion is the $\Z$-torsion in $H_2(W;\Lambda)$, which is 0.\newline
Now let us apply the universal coefficient spectral sequence
$$\Ext^p_\Lambda(H_q(W;\Lambda),\Lambda)\Longrightarrow H^{p+q}(W;\Lambda).$$
We have $H_1(W;\Lambda)=0$,  
$\Ext^p_\Lambda(H_2(W;\Lambda),\Lambda)=\Ext^p_\Lambda(\Z_-^m,\Lambda)=0$ for $p>0$,
and 
$\Ext^p_\Lambda(H_0(W;\Lambda),\Lambda)=\Ext^p_\Lambda(\Z_+,\Lambda)=0$ for $p>0$,
which implies that $$H^3(W;\Lambda)\cong \Hom_\Lambda(H_3(W;\Lambda),\Lambda).$$
Again, the right hand side is naturally a right $\Lambda$-module, and we must use the involution on $\Lambda$ to turn it into
a left $\Lambda$-module, i.e. we multiply the $T$-action with $(-1)$.
By standard theory, the map $H_3(W;\Lambda)\to H_3(W,\partial W; \Lambda)\cong \Hom_\Lambda(H_3(W;\Lambda),\Lambda)$ in the middle dimension defines the intersection product, i.e.
the map sends $x\mapsto (y\mapsto \lambda(y,x))$. 
(Observe that with our conventions, the intersection product $\lambda$ is conjugate-linear in the first and linear in the second variable.)

\bigskip

\noindent Similarly, the map $H_3(W;\Z_2)\to H_3(W,\partial W; \Z_2)\cong \Hom_{\Z_2}(H_3(W;\Z_2),\Z_2)$ 
defines the $\Z_2$-valued intersection product, i.e.
the map sends $x\mapsto (y\mapsto \lambda_{\Z_2}(y,x))$.
%The map $$H_3(W,\partial W; \Lambda)\cong \Hom_\Lambda(H_3(W;\Lambda),\Lambda)_- \to 
%H_3(W,\partial W; \Z_2)$$ factors (by naturality of the above universal coefficient spectral sequence
%for the map $\Lambda \to \Z_2$) through $\Hom_\Lambda(H_3(W;\Lambda),\Z_2)$.
 
\bigskip

\noindent We also consider the homology universal coefficient spectral sequence:
$$\Tor_p^\Lambda(H_q(W;\Lambda),\Z_2)\Rightarrow H_{p+q}(W;\Z_2).$$
This can also be interpreted as the Leray-Serre
spectral sequence for the fibration $V\to W \to \R \Pro$. 

The following diagram shows the $E^2$-term with possible differentials in low degrees.
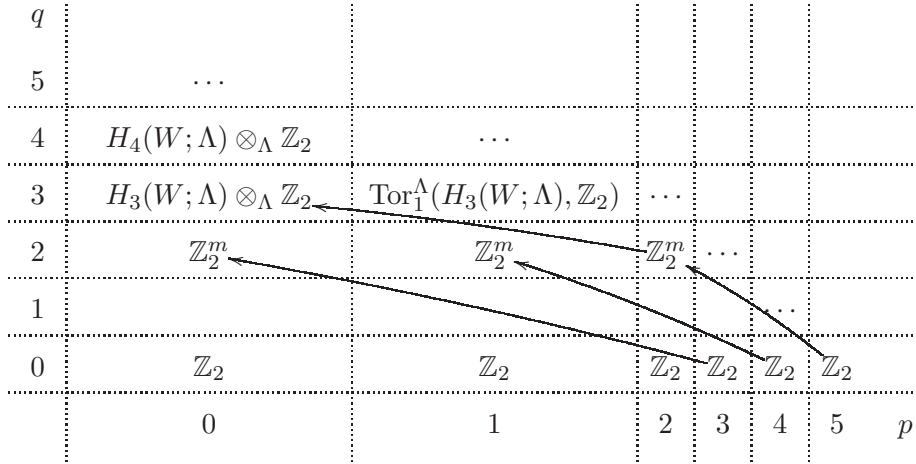
\[
 \xy <1.8pc,0pc>:<0pc,1.8pc>::
 (-9,0) *!{0},
 (-9,1) *!{1},
 (-9,2) *!{2},
 (-9,3) *!{3},
 (-9,4) *!{4},
 (-9,5) *!{5},
 (-9,6.2) *!{q},
 (-6,-1) *!{0},
 (-1,-1) *!{1}, 
 (2,-1) *!{2},
 (3,-1) *!{3},
 (4,-1) *!{4},
 (5,-1) *!{5},
 (6.2,-1) *!{p},
 (-6,0) *!{\Z_2},
 (-6,2) *!{\Z_2^m} !C="id02",
 (-6,3) *!{H_3(W;\Lambda)\otimes_\Lambda \Z_2} !C="id03",
 (-6,4) *!{H_4(W;\Lambda)\otimes_\Lambda \Z_2},
 (-1,0) *!{\Z_2},
 (-1,2) *!{\Z_2^m} !C="id12",
 (-1,3) *!{\Tor_1^\Lambda(H_3(W;\Lambda),\Z_2)},
 (2,0) *!{\Z_2},
 (2,2) *!{\Z_2^m} !C="id22",
 (3,0) *!{\Z_2} !C="id30",
 (4,0) *!{\Z_2} !C="id40",
 (5,0) *!{\Z_2} !C="id50",
 (-6,5) *!{\cdots},
 (-1,4) *!{\cdots},
 (2,3) *!{\cdots},
 (3,2) *!{\cdots},
 (4,1) *!{\cdots},
 "id30", {\ar@/_.5ex/@{->} "id02"},
 "id40", {\ar@/_.5ex/@{->} "id12"},
 "id50", {\ar@/_.5ex/@{->} "id22"},
 "id22", {\ar@/_.5ex/@{->} "id03"},
 (-9.5,-0.3) ; (6.5,-0.3) **@{.},
 (-9.5,0.7) ; (6.5,0.7) **@{.},
 (-9.5,1.7) ; (6.5,1.7) **@{.},
 (-9.5,2.7) ; (6.5,2.7) **@{.},
 (-9.5,3.7) ; (6.5,3.7) **@{.},
 (-9.5,4.7) ; (6.5,4.7) **@{.},
 (-8.5,-1.5) ; (-8.5,6.5) **@{.},
 (-3.5,-1.5) ; (-3.5,6.5) **@{.},
 (1.5,-1.5) ; (1.5,6.5) **@{.},
 (2.5,-1.5) ; (2.5,6.5) **@{.},
 (3.5,-1.5) ; (3.5,6.5) **@{.},
 (4.5,-1.5) ; (4.5,6.5) **@{.},
 
% (-.5,-.4)(0,6.1)(6,0),
 \endxy
\]

The map $H_3(W;\Lambda)\to H_3(W;\Z_2)$ can be factorized through $H_3(W;\Lambda)\otimes_\Lambda \Z_2$. The map $H_3(W;\Lambda)\to H_3(W;\Lambda)\otimes_\Lambda \Z_2$ is a surjection, and 
$H_3(W;\Lambda)\otimes_\Lambda \Z_2 \to H_3(W;\Z_2)$ can be observed at the edge 
in the universal coefficient spectral sequence. 
We can compare this to the corresponding sequence for $Q_1$. Recall that we have maps $W\to Q_m\to Q_1$ inducing isomorphisms on $\pi_1$. We will use the naturality of the universal coefficient spectral sequence.

Note that $\Tor_p^\Lambda(\Z_+,\Z_2)\cong \Tor_p^\Lambda(\Z_-,\Z_2) \cong \Z_2$ for all $p\ge 0$. Thus  
for all $p\ge 0$, we have $\Tor_p^\Lambda(H_q(Q_1;\Lambda),\Z_2)=0$ if $q$ is odd, and
$\Tor_p^\Lambda(H_q(Q_1;\Lambda),\Z_2)=\Z_2$ if $q$ is even. (The universal cover of $Q_1$ is $\C \Po^\infty$.)
We know that $H_3(Q_1;\Z_2)=0$ and that the non-zero class in $H_4(Q_1;\Z_2)$ comes from $H_4(Q_1;\Lambda)$.
So no non-zero differential in the spectral sequence for $Q_1$ ends up in $E^{0,4}$,  
and the $d_3$-differentials $E_3^{3,0} \to E_3^{2,0}$,  $E_3^{4,0} \to E_3^{2,1}$ and $E_3^{5,0} \to E_3^{2,2}$ are 
all non-zero.
By naturality, the same must be true in the spectral sequence for $W$.
It follows that $d_2:\Tor_2^\Lambda(\Z_-^m,\Z_2)=\Z_2^m \to H_3(W;\Lambda)\otimes_\Lambda \Z_2$ can have rank at most $m-1$, 
and that it is the only possibly non-zero differential ending up in $E_{0,3}$.

This implies that the map $H_3(W;\Lambda)\otimes_\Lambda \Z_2 \to H_3(W;\Z_2)$ has a cokernel of rank exactly $m-1$ (which is the sum of the ranks 
of the spaces $E_\infty^{2,1}$, $E_\infty^{1,2}$, and $E_\infty^{3,0}$ in the spectral sequence), and a kernel of rank at most $m-1$.
 
But we saw that the image of $H_3(\partial W;\Z_2)$ in $H_3(W;\Z_2)$ is the radical of the intersection 
form and has rank $m-1$.
And the map $H_3(W;\Lambda)\to H_3(W;\Z_2)$ induces a surjection onto $H_3(W;\Z_2)/rad$. 
So the images of $H_3(W;\Lambda)$ and $H_3(\partial W;\Z_2)$ in $H_3(W;\Z_2)$ must be complements.
This implies that the image of $H_3(W;\Lambda)$ in $H_3(W;\Z_2)$
injects into $H_3(W,\partial W;\Z_2)$.

But this also implies that all $\Lambda$-torsion in $H_3(W;\Lambda)$ maps to $0$ in $H_3(W, \partial W;\Z_2)$,
respectively in $H_3(W;\Z_2)$: Since the $\Z_2$-intersection form is 0 on a complement of the image of $H_3(W;\Lambda)$,
it suffices to show that if $x\in H_3(W;\Lambda)$ is torsion, then $\epsilon(\lambda(x,y))=0$ for all $y \in H_3(W;\Lambda)$.
We have three cases: If $x$ is $\Z$-torsion, then $\lambda(x,y)=0$ for all $y$.
If $(1+T)x=0$, then $T\lambda(x,y)=\lambda(-Tx,y)=\lambda(x,y)$, so $\lambda(x,y)$ is a multiple of $(1+T)$ and  $\epsilon(\lambda(x,y))=0$ for all $y$. 
And similarly if $(1-T)x=0$, then $T\lambda(x,y)=\lambda(-Tx,y)=\lambda(-x,y)=-\lambda(x,y)$, so $\lambda(x,y)$ is a multiple of $(1-T)$ and  $\epsilon(\lambda(x,y))=0$ for all $y$.

\bigskip

Let $H_1(M)\cong S\oplus \Z^{m-1}$, where $S$ is odd torsion. Then $H_2(M)\cong \Z^{m-1}$.
Now the structure theorem for finitely generated $\Z$-free $\Lambda$-modules
\cite{CR} says that every such module is a direct sum of summands isomorphic to 
$\Lambda$, $\Z_+$ or $\Z_-$.
We want to apply this to the exact sequence 
$$0\to S_-\oplus \Z_-^{m-1} \to H_3(W;\Lambda) \to H_3(W, \partial W;\Lambda) \to \Z_+^{m-1}\to 0.$$
Recall that $H_3(W, \partial W;\Lambda)$ is free over $\Z$, which implies that over $\Z$, 
we get two split short exact sequences $0\to S\oplus \Z^{m-1} \to H_3(W;\Lambda) \to K\to 0$ and
$0\to K \to H_3(W, \partial W;\Lambda) \to \Z_+^{m-1}\to 0$, where $K$ is the image of
the map  $H_3(W;\Lambda) \to H_3(W, \partial W;\Lambda)$.

By applying the structure theorem, we see that the modules $H_3(W;\Lambda)/\{\Z -torsion\}$ and $H_3(W, \partial W;\Lambda)$ have the form
$H_3(W;\Lambda)/\{\Z -torsion\} \cong \Z_-^a\oplus \Z_+^b\oplus \Lambda^c$,
and $H_3(W,\partial W;\Lambda)\cong \Z_+^a\oplus \Z_-^b\oplus \Lambda^c$.
We have also used that $H_3(W,\partial W;\Lambda)\cong \Hom_\Lambda(H_3(W;\Lambda)/\{\Z -torsion\},\Lambda)$ with $T$-action multiplied by $-1$.
Note that only the numbers $a,b,c$ are uniquely determined; 
the decomposition itself into cyclic summands is NOT unique. 

Choose such a decomposition for $H_3(W;\Lambda)/\{\Z -torsion\}$, and use the dual decomposition for $\Hom_\Lambda(H_3(W;\Lambda),\Lambda)$. 
By tensoring the exact sequence over $\Lambda$ with $\Q[\Z_2]\cong \Q_-\oplus \Q_+$, we get
$$0\to \Q_-^{m-1} \to \Q_-^{a+c}\oplus \Q_+^{b+c} \to \Q_-^{b+c}\oplus \Q_+^{a+c}\to \Q_+^{m-1}\to 0.$$
It follows that $a=b+(m-1)$.
So $H_3(W;\Lambda)/\{\Z -torsion\} \cong \Z_-^{b+(m-1)}\oplus \Z_+^b\oplus \Lambda^c$.

Now we really use the assumption that $S$ consists only of odd torsion elements: Then $Ext_\Lambda^1(\Z_\pm, S_-)=0$, and we obtain
$H_3(W;\Lambda) \cong S_- \oplus \Z_-^{b+(m-1)}\oplus \Z_+^b\oplus \Lambda^c$.
Then 
\begin{eqnarray*}
rk(H_3(W;\Lambda)\otimes_\Lambda \Z_2) &=& rk(S_- \oplus \Z_-^{b+(m-1)}\oplus \Z_+^b\oplus \Lambda^c)\otimes_\Lambda \Z_2 \\
& = & 0 + b + (m-1) + b + c = 2b + c + (m-1).
\end{eqnarray*}
We have seen that the map $H_3(W;\Lambda)\to H_3(W;\Z_2)$ is 0 on the submodule $S_-\oplus \Z_-^{b+(m-1)}\oplus \Z_+^b$ of 
$H_3(W;\Lambda)$ and that the map $H_3(W;\Lambda)\otimes_\Lambda \Z_2\to H_3(W;\Z_2)$ has kernel of rank at most $m-1$.
It follows that $b=0$, and since the image has even rank, $c$ must be even, say $c=2k$.
So our exact sequence becomes
$$0\to S_- \oplus \Z_-^{m-1} \to S_- \oplus \Z_-^{m-1} \oplus \Lambda^{2k} \to \Z_+^{m-1} \oplus \Lambda^{2k} \to \Z_+^{m-1}\to 0.$$

The $\Z$-torsion $S$ is mapped isomorphically to itself under the first map, and the intersection form is $0$ if one of the
arguments is in $S$. So let us again look at the sequence modulo the $\Z$-torsion $S$: 
$$0\to  \Z_-^{m-1} \stackrel{(i_1,i_2)}\longrightarrow \Z_-^{m-1} \oplus \Lambda^{2k} \to \Z_+^{m-1} \oplus \Lambda^{2k} \to \Z_+^{m-1}\to 0.$$
On the submodule $\Lambda^{2k} \subseteq H_3(W;\Lambda)$, we know that if we tensor the quadratic form with $\Z_2$,
it becomes hyperbolic. So we can choose basis vectors $e_i,f_i$ for $\Lambda^{2k}$ such that 
$\epsilon(\lambda(e_i,e_j))=\epsilon(\lambda(f_i,f_j))=0$ and 
$\epsilon(\lambda(e_i,f_j))=\delta_{ij}$.

\noindent What can we say about the matrix $A$ of the intersection form $\lambda$ with respect to the $e_i,f_j$?
\begin{itemize}
\item
$A$ is a skew-hermitian matrix.
\item
$A$ is the matrix for the component $\Lambda^{2k} \to \Lambda^{2k}$ of the middle map in the exact sequence.
\item
$\epsilon(A)=\left( \begin{array}{cc} 0 & I\\ I & 0 \end{array} \right)$. This means in particular that 
$A$ has a determinant which is a nonzerodivisor, which implies that the map $\Lambda^{2k} \to \Lambda^{2k}$ is injective.
\end{itemize}
So we see that the middle map in the spectral sequence is injective on $\Lambda^{2k}$.
Let us denote its image by $L \subseteq \Z_+^{m-1} \oplus \Lambda^{2k}$. We have $L\cong \Lambda^{2k}$ as $\Lambda$-modules.
We also get that $i_1$ is injective, thus it has a finite cokernel, which we denote by $G_-$.

\bigskip

Our next claim is that $G_-$ does not have non-trivial $2$-torsion.
Let $v\in \Z_-^{m-1}$ be a representative for a $2$-torsion element in the cokernel, and
suppose that $x=(2v,-w)\in \Z_-^{m-1} \oplus \Lambda^{2k}$ is in the image of $(i_1,i_2)$. This is equivalent to 
$\lambda(2v,y)=\lambda(w,y)$ for all $y\in H_3(W;\Lambda)$.
It suffices to show that $w$ is divisible by 2, because then $(v,-w/2)$ is in the image of $(i_1,i_2)$, and 
$v$ represents $0\in G_-$.

Note that $\lambda(2v,y)$ is divisible by $2(1+T)$ for all $y$ and that $Tw=-w$, since it is in the image of a map from $\Z_-$. Thus we can
write $$w=\sum \alpha_i (1-T) e_i + \beta_i (1-T) f_i,$$ where $\alpha_i,\beta_i\in \Z$.
Set $y=f_j$. Then $\lambda(2v,y)=\lambda(w,y)$ implies that 
$$\sum \alpha_i (1+T) \lambda(e_i,f_j) + \beta_i (1+T) \lambda(f_i,f_j)$$ is divisible by $2(1+T)$.
Now for $i\ne j$, we have $\epsilon (\lambda(e_i,f_j))=0$ , and this implies that $(1+T) \lambda(e_i,f_j)$ is divisible by $2(1+T)$.
Similarly, $(1+T) \lambda(f_i,f_j)$ is divisible by $2(1+T)$ for all $i$.
But $(1+T) \lambda(e_j,f_j)$ is not divisible by 2, and so $\alpha_j$ is divisible by $2$.
Repeating the argument with $y=e_j$ shows that $\beta_j$ is divisible by 2.
Thus we have shown that $w$ is divisible by 2.

\bigskip

By a diagram chase, we see that the cokernel $G_-$ of $i_1$ is isomorphic to $K/L$. (Recall that $K$ is the image of the middle map, $L$ is the image of the middle map restricted to the free submodule.) 
It follows that $K$ is an extension of the finite $\Lambda$-module $G_-$  by the free $\Lambda$-module $L\cong \Lambda^{2k}$. We must have $K\cong \Z^{4k}$ as $\Z$-module, since $K$ is a submodule of $H_3(W,\partial W; \Lambda)$. Thus $K$ is a direct sum of summands isomorphic to $\Lambda$, $\Z_+$ or $\Z_-$ as $\Lambda$-module.

If we tensor the exact sequence $$0\to \Lambda^{2k}\to K \to G_-\to 0$$ over $\Lambda$ with $\Z_2$, we get the
(right) exact sequence $$\Z_2^{2k}\to K\otimes_\Lambda \Z_2 \to 0\to 0.$$ (We have $G_-\otimes_\Lambda \Z_2=0$
since $G_-$ has trivial $2$-torsion.) But since $K\otimes_\Lambda \Z_2$ has rank $\le 2k$, $K$ must be a free $\Lambda$-module, i.e. $K\cong \Lambda^{2k}$.

But since $K$ is free, this means that the exact sequence 
$$0\to  \Z_-^{m-1} \stackrel{(i_1,i_2)}\longrightarrow \Z_-^{m-1} \oplus \Lambda^{2k} \to K\to 0$$ 
splits. Thus there is a new decomposition $H_3(W;\Lambda)/S \cong \Z_-^{m-1}\oplus \Lambda^{2k}$
such that the first summand $\Z_-^{m-1}$ is exactly the kernel for the map to $H_3(W,\partial W;\Lambda)$.
Choose basis vectors $a_j$ for this submodule. We may choose basis vectors for the second summand with the same properties as the $e_i,f_i$. By abuse of notation, and since we do not need the old generators any more, we denote these new basis vectors by $e_i,f_i$ again.

\bigskip

Then $\lambda(a_j,y)=0$ for all $y\in H_3(W;\Lambda)$. 
The generators $a_j,e_i,f_i$ for $$H_3(W;\Lambda)/S \cong \Z_-^{m-1}\oplus \Lambda^{2k}$$ can be used to give generators
for $H_3(W,\partial W;\Lambda)\cong \Hom_\Lambda(H_3(W;\Lambda),\Lambda)$.
For each $j$ we have a unique such map sending $a_j$ to $(1-T)$ and all other generators to $0$.
For each $i$ we have a unique such map sending $e_i$ to $1$ and all other generators to $0$
and a unique such map sending $f_i$ to $1$ and all other generators to $0$.
These maps form generators for $H_3(W,\partial W;\Lambda)\cong \Z_+^{m-1} \oplus \Lambda^{2k}$, i.e. our distinguished
decomposition of $H_3(W;\Lambda)$ induces a distinguished decomposition of $H_3(W,\partial W;\Lambda)$.

Now we see that $\Lambda^{2k} \subseteq H_3(W;\Lambda)$ has to map into $\Lambda^{2k} \subseteq H_3(W,\partial W;\Lambda)$
since $\lambda(a_j,y)=0$ for all $y$ also implies $\lambda(y,a_j)=0$ for all $y$.
The map $\Lambda^{2k} \to \Lambda^{2k}$ must have zero cokernel, so it is bijective.

\bigskip

So $H_3(W;\Lambda)/rad$ is a free $\Lambda$-module generated by $e_i,f_i$ and the induced quadratic form is 
non-degenerate.
But this implies we finally can use the standard theory of non-degenerate quadratic forms. The only 
difference is that our quadratic refinement $\tilde{\mu}$ has values in $\Lambda / \Z \cdot 1$, whereas
Wall's quadratic refinement takes values in 
$\Lambda / \langle x + \bar{x} \mid x\in \Lambda \rangle = \Lambda / \Z \cdot 2$.
Wall proved that $L_6(\Lambda)\cong \Z_2$, generated by the Arf invariant \cite{SCM}. 
This implies that $\tilde{L_6}(\Lambda)$ is the trivial group: there is a lift of $\tilde{\mu}$
from $\Lambda / \Z \cdot 1$ to $\Lambda / \Z \cdot 2$ with zero Arf invariant.

\bigskip

Thus, after stabilizing $W$ by connected sum with copies of $S^3\times S^3$, we may assume that the restriction of $\lambda$ 
to the free submodule generated by $e_i,f_i$ is hyperbolic.
We may also assume that the $e_i$ are generators for a Lagrangian. Since the Hurewicz map $\pi_3(W)\to H_3(W;\Lambda)$
is surjective, we obtain disjointly embedded spheres in $W$ with trivial normal bundle, more precisely with a unique trivialization
such that the normal $B$-structure extends over the surgery cobordism.

\bigskip

Now the diagram on p. 73 of Ranicki's book \cite{Ran} shows that the effect of doing surgery on the $e_i$ 
on the homology of $W$ with coefficients in $\Lambda$ is exactly killing the free submodule  of $H_3(W;\Lambda)$ generated by the $e_i,f_i$. 
In particular, the resulting manifold $W_0$ still has $2$-type $B$.
Thus the Mayer-Vietoris sequence for $X=M\times D^3 \cup V_0$ contains
$$H_2(M\times S^2) \to H_2(M\times D^3) \oplus H_2(V_0) \to H_2(X) \to 0.$$ 
We deduce that $H_2(X) \cong \Coker(H_2(S_2)\to H_2(V_0))$.
But in \cite{OlbThesis} we have seen that if $W_0$ has normal 2-type $B$, then this cokernel is isomorphic to $\Z_-^{m-1}$.

So we have constructed a conjugation on a simply connected spin manifold $X$ with fixed point set $M$.
But since $H_1(X)$ and $H_2(X)$ are free over $\Z$, Poincar\'e duality implies that all homology of $X$ is free over $\Z$.
And since we have a degree-halving $\Z_2$-homology isomorphism from $X$ to $M$, we see that the free homology of $X$ is concentrated in even degrees.

\section{Analysis of the result of the construction}
\subsection{The $\Z_2$-cohomology and homology of $Q_m$}

The integral cohomology of $K(\Z^m,2)$ is a polynomial ring on the standard
classes $v_1,\dots, v_m\in H^2(K(\Z^m,2))$ which correspond to the dual basis of the standard basis of $\Z^m$.
We get dual generators $e_i=v_i^*$ for the second integral homology corresponding to the standard basis of $\Z^m$.
More generally we get a basis $e_{i_1\dots i_k}=(v_{i_1}\dots v_{i_k})^*$ for the integral homology in degree $2k$.
Let us abuse notation and denote by $v_i$ also the corresponding class in
$\Z_2$-cohomology, and by $e_{i_1\dots i_k}$ also the corresponding class 
in $\Z_2$-homology.

In \cite{OlbThesis} it was proved that 
$$H^*(Q_m;\Z_2)\cong \Z_2[q,t,x_1,\dots x_{m-1}]/t^3\quad \text{with }deg(q)=4, deg(t)=1, deg(x_i)=2.$$
Here we have to be more precise about the isomorphism.
We denote the maps 
\begin{eqnarray*}
 Q_m = (\C \Pro \times \dots \times \C \Pro)/\tau^m & \to & Q_n = (\C \Pro \times \dots \times \C \Pro)/\tau^n \\
 \left[ a_1,\dots, a_m \right] & \mapsto & \left[ a_{i_1},\dots, a_{i_n} \right]
\end{eqnarray*}
by $pr_{i_1,\dots i_n}$.
There is a fibration $\R \Po^2\to Q_1 \to \HB \Pro$ whose Serre spectral sequence collapses at the $E_2$-term. 
We get generators $t\in H^1(Q_1;\Z_2)$ and $q\in H^4(Q_1;\Z_2)$. 
Using the fibration $(\C \Pro)^{m-1}\to Q_m \stackrel{pr_m}{\to} Q_1$ we transport the classes
$t$ and $q$ to $H^*(Q_1;\Z_2)$. For $i=1,\dots m-1$ we denote by $x_i\in H^2(Q_m)$ the class which 
maps to $v_i\in H^2(K(\Z^{m-1},2))$, and to $0$ under the diagonal map $\Delta=pr_{1,\dots,1}: Q_1\to Q_m$.   

To compute the maps $\pi:K(\Z^m,2)\to Q_m$ and $pr_i: Q_m\to Q_1$ in $\Z_2$-cohomology for $i<m$, we first look at the case $m=2$. 
We consider the commutative diagram whose rows are fibrations:
$$ 
\xymatrix{
S^2\times S^2
\ar[d]
\ar[r]
& \C \Pro \times \C \Pro
\ar[d]
\ar[r]
& \HB \Pro \times \HB \Pro
\ar[d]^{id}
\\
(S^2\times S^2)/(-1,-1)
\ar[r]
& Q_2
\ar[d]^{pr_1,pr_2}
\ar[r]
& \HB \Pro \times \HB \Pro
\\
& Q_1\times Q_1
\ar[ur]
}
$$
The elements $q_1,q_2\in H^4(\HB\Pro\times \HB\Pro;\Z_2)$ (coming from the two factors)
pull back to $v_1^2,v_2^2\in H^4(\C\Pro\times \C\Pro;\Z_2)$, and to $0\in  H^4((S^2\times S^2)/(-1,-1);\Z_2)$.
Thus under $\pi:K(\Z^2,2)\to Q_2$, the element $q$ maps to $v_2^2$, the element $t$ maps to 0, and the element $x$ maps to $v_1+v_2$.
And under $pr_1:Q_2\to Q_1$, the element $t$ maps to $t$, and $q$ maps to $q+x_1^2$ 
(there is no term $t^2x_1$ since this would map nontrivially to $(S^2\times S^2)/(-1,-1)$).
 
Now for the general case we consider the diagram (where $i<m$):
$$ 
\xymatrix{
& K(\Z^m,2)
\ar[d]
\ar[r]
& K(\Z^2,2)
\ar[d]
\\
Q_1
\ar[r]^{\Delta}
& Q_m
\ar[r]^{pr_{i,m}}
\ar[dr]_{(pr_i, pr_m)}
& Q_2
\ar[d]
\\
& & Q_1\times Q_1
}
$$
Here $q_1,q_2\in H^4(Q_1\times Q_1;\Z_2)$ map to $q+x_1^2,q\in H^4(Q_2;\Z_2)$
and to $v_i^2, v_m^2\in H^4(K(\Z^m,2);\Z_2)$. 
Hence $x_1\in H^2(Q_2;\Z_2)$ maps to $x_i\in H^2(Q_m;\Z_2)$ under $pr_{i,m}$.
Under $pr_i:Q_m\to Q_1$, we have $t\mapsto t$, $q\mapsto q+x_i^2$,
and under $\pi:K(\Z^m,2)\to Q_m$, we get $t\mapsto 0, x_i\mapsto v_i+v_m, q\mapsto v_m^2$.

Similarly, we compute the effect of $pr_{i,j}: Q_m \to Q_2$ for $i<j<m$. We obtain
$q\mapsto q+x_j^2$, $x_1\mapsto x_i+x_j$, $t\mapsto t$.
And for $pr_{ijk}:Q_m\to Q_3$ with $i<j<k$ we get 
$q\mapsto q+x_k^2$, $x_1\mapsto x_i+x_k$, $x_2\mapsto x_j+x_k$, $t\mapsto t$ in the case $k<m$ 
and $q\mapsto q$, $x_1\mapsto x_i$, $x_2\mapsto x_j$, $t\mapsto t$ in the case $k=m$. 

\bigskip

We will use monomials in the generators as standard bases for the 
$\Z_2$-cohomology groups and use the dual bases for the corresponding
$\Z_2$-homology groups. For example $(qx_1)^*$ will denote the element of 
the sixth $\Z_2$-homology group which pairs to $1$ with $qx_1$ and to $0$ 
with all other monomials in the generators.

\bigskip

The transfer map in $\Z_2$-homology $tr:H_*(Q_m;\Z_2)\to H_*(K(\Z^m,2);\Z_2)$ is induced by the change of coefficients (modules over the group ring
of the fundamental group of $Q_m$) given by $\Z_2\to \Z_2[\Z_2]$, $1\mapsto 1+T$. 
We use the short exact sequence $\Z_2\to \Z_2[\Z_2] \to \Z_2$ which induces a long exact sequence in homology. The fact that 
$\Z_2[\Z_2] \to \Z_2$ induces $\pi_*:H_*(K(\Z^m,2);\Z_2)\to H_*(Q_m;\Z_2)$ shows that $Im(tr)=Ker(\pi_*)$.
From the fact that the involution $\tau^m$ on $K(\Z^m,2)$ is trivial on $\Z_2$-homology it follows that $Im(\pi_*)\subseteq Ker(tr)$. 
(The composition $\Z_2[\Z_2] \to \Z_2 \to \Z_2[\Z_2]$ is multiplication by $1+T$, which induces $id+(\tau^m)_*$.)
But since the ranks of $H_*(Q_m;\Z_2)$ equal those of $H_*(K(\Z^m,2);\Z_2)$ in even degrees, we get equality 
$Im(\pi_*)= Ker(tr)$ in these degrees.

We first look at $tr:H_2(Q_m;\Z_2)\to H_2(K(\Z^m,2);\Z_2)$. 
From the above computation of $\pi^*$ it follows that $tr((x_i)^*)=0$ for all $i=1, \dots ,m-1$
and $tr((t^2)^*)=\sum_{j=1}^m e_j$.

The kernel of $\pi_*:H_4(K(\Z^m,2);\Z_2)\to H_4(Q_m;\Z_2)$ is generated by
$\sum_{j, j\ne i} e_{ij}, i=1,\dots, m-1$ (just dualize $\pi^*$),
and the image of $\pi_*$ has generators $q^*, (x_ix_j)^*$ for $i\le j$. 
To compute the transfer map in degree 4, we first deduce that for $m=2$,
the remaining generator $(t^2x_1)^*$ must be mapped to $e_{12}$.
Then we use naturality with respect to the projections $pr_{j,k}$ to see that
for general $m$, the transfer map is $(t^2x_i)^*\mapsto \sum_{j, j\ne i} e_{ij}$.

In degree 6, the transfer has image generated by
$$
\sum_{i\le j\le k} e_{ijk},\quad
\sum_{j=1}^m e_{iij}\text{ for }i=1,\dots, m-1, \quad
\sum_{\stackrel{k\ne i,j}{k=1}}^m e_{ijk}\text{ for $i<j<m$},
$$
and kernel generated by
$$ (x_iq)^*, (x_ix_jx_k)^* \text{ for $i\le j\le k$},$$ 
which follows again from our earlier computation of $\pi^*$.
Again we first compute the cases for small $m\le 3$, and then
use the naturality of the transfer with respect to projections to $Q_1$, $Q_2$ and $Q_3$
to see that the transfer is 
\begin{gather*}
(t^2q)^*\mapsto \sum_{i\le j\le k} e_{ijk},
(t^2x_i^2)^*\mapsto \sum_{j=1}^m e_{iij}\text{ for }i=1,\dots, m-1,  \\
(t^2x_ix_j)^*\mapsto \sum_{\stackrel{k\ne i,j}{k=1}}^m e_{ijk}\text{ for $i<j<m$}.
\end{gather*}

\subsection{The transfer map $\Omega_6^B\to \Omega_6^{\tilde{B}}$}

The projection map $\pi_*:\Omega_6^{\tilde{B}}\to \Omega_6^B$ is given 
by $[M\to \tilde{B}]\mapsto [M\to \tilde{B}\to B]$. The composition
$tr\circ \pi_*$ is equal to $1+\tau_*$, where $\tau$ denotes the nontrivial 
deck transformation $\tilde{B}\to {\tilde{B}}$. This is multiplication with $-1$ on $H^2(\tilde{B})$, 
but it also changes the orientation of $M$, because the orientation is obtained by composition with $\tilde{B}\to BSO$,
and the nontrivial deck transformation of the double cover $BSO\to BO$ reverses orientation of a bundle. 
By the computation of $\Omega_6^{\tilde{B}}$ we see that $\tau_*$ is the identity. 
So $tr\circ \pi_*$ is multiplication by 2.
In particular it follows that the cokernel of $tr$ is a finite group consisting
of elements of order (dividing) 2.

We can compute the transfer map using the Atiyah-Hirzebruch spectral sequence.
$$ 
\xymatrix{
E^2_{p,q}\cong H_p(Q_m;\underline{\Omega_q^{Spin}}) 
\ar[d]
\ar[r]
& \Omega_6^B
\ar[d]
\\
\tilde{E}^2_{p,q}\cong H_p(K(\Z^m,2);{\Omega_q^{Spin}}) 
\ar[r]
& \Omega_6^{\tilde{B}}
}
$$

For the facts needed and not proved here see Zubr \cite{Zub} (for the computation of
$\Omega_6^{\tilde{B}}$) and \cite{OlbThesis} (for the computation of 
$\Omega_6^{B}$).

For the terms on the second pages of the spectral sequences
we use the transfer maps in homology induced by the short exact coefficient
sequences $\Z_-\to \Z[\Z_2] \to \Z$ and $\Z_2\to \Z_2[\Z_2] \to \Z_2$.
We get:
$$E^2_{0,6}=\tilde{E}^2_{0,6}=E^2_{1,5}=\tilde{E}^2_{1,5}=0.$$
The map $$E^2_{2,4}\cong \Z^m \to \tilde{E}^2_{2,4}\cong \Z^m$$
is injective with cokernel $H_2(Q_m)\cong \Z_2^{m-1}$. More precisely its 
image is the kernel of the map $H_2(K(\Z^m,2))\cong \Z^m\to H_2(Q_m)\cong
\Z_2^{m-1}$, which is generated by $2 e_i, i=1, \dots, m-1$ and 
$\sum_{j=1}^m e_j$. (Compare with the $\Z_2$-homology transfer.)

We have $$E^2_{3,3}\cong 0 \to \tilde{E}^2_{3,3}\cong 0.$$
We do already know the $\Z_2$-homology transfer map  
$$E^2_{4,2}\cong \Z_2^{{m+1}\choose 2} \to 
\tilde{E}^2_{4,2}\cong \Z_2^{{m+1}\choose 2}.$$
We have
$$E^2_{5,1}\cong \Z_2^{{m\choose 2}+1} \to 
\tilde{E}^2_{5,1}\cong 0.$$
Finally we see that
$$E^2_{6,0}\cong \Z^{{m+2}\choose 3} \to 
\tilde{E}^2_{6,0}\cong \Z^{{m+2}\choose 3}$$
is injective with cokernel isomorphic to $H_6(Q_m)\cong \Z_2^{{{m+1}\choose 3}+m-1}$.
More precisely the image in 
$E^2_{6,0}\cong H_6(K(\Z^m,2))$ is the kernel of the composition
$H_6(K(\Z^m,2))\to H_6(K(\Z^m,2);\Z_2)\to H_6(Q_m;\Z_2)$. 
We obtain generators 
\begin{gather*}
\sum_{j=1}^m e_{iij}\text{ for }i=1,\dots, m-1, \quad \sum_{\stackrel{k\ne i,j}{k=1}}^m e_{ijk}\text{ for $i<j<m$},\\
\sum_{i\le j\le k} e_{ijk},\quad 2 e_{ijk} \text{ for }i\le j\le k < m,\quad  2e_{imm}\text{ for $i<m$}.
\end{gather*} 

For the third pages of the spectral sequences we get the same as before
in the $(2,4)$-terms. We get the zero map $E^3_{4,2}\cong \Z_2^{m} \to 
\tilde{E}^3_{4,2}\cong \Z_2^{m}$, since the image of the map on the second page is
also in the image of the $\tilde{d}_2$-differential. 
We get a zero map $E^3_{5,1}\cong \Z_2 \to \tilde{E}^3_{5,1}\cong 0$. 

Finally for the
$(6,0)$-terms, we have to restrict both domain and codomain of the second page
to subgroups with quotient $\Z_2^{m\choose 2}$:

The $\tilde{d}_2$-differential is the composition of reduction modulo 2
and the dual of $Sq^2$. It sends $e_{ijj}$ for $i\ne j$ to $e_{ij}$
and all other generators to 0.
Thus we can consider $\tilde{E}^3_{6,0}$ as the subset of 
$\tilde{E}^2_{6,0}$ with generators
$e_{ijk}$ for $i<j<k$, $e_{iii}$ for all $i$,
$e_{iij}+e_{ijj}$ and $2e_{iij}$ for all $i<j$.
We saw that $E^2_{6,0}$ can also be considered as a subset of $\tilde{E}^2_{6,0}$.
We have to compute the image of the generators under $\tilde{d}_2$.

The $d_2$-differential is the composition of reduction modulo 2
and the dual of $Sq^2+tSq^1$. Since reduction modulo 2 is surjective, its image
is the image of the dual of $Sq^2+tSq^1$. Its generators are 
$(x_ix_j)^*$ for $i<j$ and $(t^2x_i)^*+(x_i^2)^*$ for all $i$.

$$ 
\xymatrix{
E^3_{6,0} \cong \Z^{{m+2}\choose 3}
\ar[d]
\ar[r]
& E^2_{6,0} \cong\Z^{{m+2}\choose 3}
\ar[r]
\ar[d]
& \langle (x_ix_j)^*, (t^2x_i)^*+(x_i^2)^* \rangle \cong \Z_2^{{m}\choose 2}
\ar[d]
\\
\tilde{E}^3_{6,0} \cong\Z^{{m+2}\choose 3}
\ar[r]
& \tilde{E}^2_{6,0}\cong \langle e_{ijk} \mid i\le j\le k \rangle \cong\Z^{{m+2}\choose 3}
\ar[r]
& \langle e_{ij} \rangle \cong \Z_2^{{m}\choose 2}
}
$$

We begin with the case $m=1$. Here $E^2_{6,0}$ is generated by $e_{111}$, and this generator maps to 0.
In the case $m=2$ we use the diagonal map and the case $m=1$ to compute $d_2(e_{111}+e_{112}+e_{122}+e_{222})=0$,
we use the mod 2 transfer computations to see that $e_{111}+e_{112}\mapsto (t^2x_1)^*+(x_1^2)^*$.
For the remaining generators we use the coefficient changes $\Z[\Z_2]\to \Z_-$ and $\Z_-\to\Z[\Z_2]$ which 
induce projection and transfer maps respectively. Their composition $\Z[\Z_2]\to \Z[\Z_2]$ is multiplication with 
$1-T$. On $H_6(K(\Z^m,2))$ this is multiplication with $2$. The other generators are in the image of the projection map
$H_6(K(\Z^m,2))\to H_6(Q_m;\Z_-)$. This allows to compute that $d_2(2e_{111})=d_2(2e_{122})=0$.

For $m=3$, we compute
$e_{111}+e_{112}+e_{113}+e_{122}+e_{123}+e_{133}+e_{222}+e_{223}+e_{233}+e_{333}\mapsto 0$
using the diagonal map,
$d_2(2e_{111})=d_2(2e_{222})=0$ and $d_2(2e_{112})=d_2(2e_{122})=d_2(2e_{133})=d_2(2e_{233})=(x_1x_2)^*$
using the projection map.
Finally we use the maps $pr_{1,1,2}:Q_2\to Q_3$, $pr_{1,2,2}:Q_2\to Q_3$ and $pr_{2,1,2}:Q_2\to Q_3$ and the case $m=2$ 
to compute $e_{111}+e_{112}+e_{113}\mapsto (t^2x_1)^*+(x_1^2)^*$, $e_{221}+e_{222}+e_{223}\mapsto (t^2x_2)^*+(x_2^2)^*+\mu (x_1x_2)^*$,
and $e_{123}\mapsto (x_1x_2)^*$. 

For the general case, we use the $\Z_2$-homology transfer computation to conclude that
$\sum_{j=1}^m e_{iij} \mapsto (t^2x_i)^*+(x_i^2)^* + \sum_{j<k}\mu_{ijk}(x_jx_k)^*$ for certain $\mu_{ijk}\in\Z_2$, 
and that all other generators map to a linear combination of generators of the form $(x_ix_j)^*$.
Using the projection map, we compute that $2e_{ijk}\mapsto 0$ in the cases $i<j<k<m$ and $i=j=k<m$, 
$2e_{ijk}\mapsto (x_ix_j)^*$ for $i<j=k<m$ and $2e_{ijk}\mapsto (x_jx_k)^*$ for $i=j<k<m$.
Also we obtain $2e_{imm}\mapsto \sum_{j\ne i}(x_ix_j)^*$.
Finally we compare with the case $m=3$ to see that
$\sum_{\stackrel{k\ne i,j}{k=1}}^m e_{ijk}\mapsto (x_ix_j)^*$ and $\sum_{i\le j\le k} e_{ijk}\mapsto 0$.

This finally allows us to identify generators for the kernel $E^3_{6,0}$:
\begin{gather*}
\sum_{i\le j\le k} e_{ijk}, \quad 
2e_{iii}\text{ for }i<m, \quad 
\sum_{j=1}^m 2e_{iij}\text{ for }i<m, \quad
2e_{imm}+\sum_{k\ne i, k\ne m} e_{imk}\text{ for }i<m, \\
2e_{ijj}+\sum_{k\ne i, k\ne j} e_{ijk}\text{ for }i\ne j\text{ and }i,j<m, \quad 
2e_{ijk}\text{ for }i<j<k.
\end{gather*}

This means that the classes $F\in H_6(K(\Z^m,2))$ belonging to $E^3_{6,0}$ are exactly those which satisfy:
\begin{eqnarray*}
\langle v_iv_iv_j, F \rangle & = & \langle v_kv_kv_l, F \rangle \ (mod\ 2) \quad \text{ for all }i,j,k,l\in\{1,\dots m\} \\
\langle v_iv_jv_k, F \rangle & = & \frac12 ( \langle v_iv_iv_j, F \rangle + \langle v_iv_jv_j, F \rangle + \langle v_iv_iv_k, F \rangle \\
& & \quad + \langle v_iv_kv_k, F \rangle + \langle v_jv_jv_k, F \rangle + \langle v_jv_kv_k, F \rangle )\ (mod\ 2) \quad \text{ for all }i<j<k.
\end{eqnarray*}

\bigskip

There are almost no more higher differentials involved.
The only term which could possibly be different on the ${\infty}$-page is 
$E^3_{4,2}$ since there is room for a $d_3$ from $E^3_{7,0}$.
Equipped with this information about the ${\infty}$-page, we compute 
the extensions. We denote the filtrations of $\Omega_6^{\tilde{B}}$ and  
$\Omega_6^{B}$ by $\tilde{F}_j$ and $F_j$ so that 
$\tilde{F}_j/\tilde{F}_{j-1}\cong \tilde{E}^{\infty}_{j,6-j}$
and $F_j/F_{j-1}\cong E^{\infty}_{j,6-j}$.

Actually we know from \cite{OlbThesis} that 
$$\Omega_6^B\cong \Z^{{{m+2}\choose 3 }+m}\oplus \Z_4$$ 
and from Zubr that $$\Omega_6^{\tilde{B}}\cong \Z^{{{m+2}\choose 3 }+m}.$$
We have a commutative diagram with short exact sequences:
$$ 
\xymatrix{
F_5 
\ar[d]
\ar[r]
& \Omega_6^B
\ar[r]
\ar[d]
& E^3_{6,0}
\ar[d]
\\
\tilde{F}_5 
\ar[r]
& \Omega_6^{\tilde{B}}
\ar[r]
& \tilde{E}^3_{6,0}
}
$$
We see that $F_5 \cong \Z^m\oplus \Z_4$ and $\tilde{F}_5 \cong \Z^m$.
We apply the snake lemma, using that the middle vertical map has cokernel a 
$\Z_2$-vector space, and the right vertical map is injective.
We get exact sequences $0\to \Z_4\to F_5 \to \tilde{F}_5$, 
$0\to \Z_4\to \Omega_6^B\to \Omega_6^{\tilde{B}}$, 
and for each of these sequences the last
map has cokernel a finite-dimensional $\Z_2$-vector space. In particular the
dimension of this vector space for the first exact sequence is at most $m$.

Now we consider the commutative diagram with short exact sequences:
$$ 
\xymatrix{
F_2 \cong \Z^m 
\ar[d]
\ar[r]
& F_5 \cong \Z^m\oplus \Z_4
\ar[r]
\ar[d]
& F_5/F_2
\ar[d]
\\
\tilde{F}_2 \cong \Z^m
\ar[r]
& \tilde{F}_5 \cong \Z^m
\ar[r]
& \tilde{F}_5/\tilde{F}_2 \cong \Z_2^{m}
}
$$
We also have the short exact sequence $E^{\infty}_{4,2}\to F_5/F_2 \to \Z_2$.
By applying the snake lemma again, we see that $F_5/F_2$ contains an element of
order $4$, so it is isomorphic to $\Z_4\oplus \Z_2^k$ for some $k\le m-1$. We
may also assume that the generator of the $\Z_4$-summand is in the kernel of 
the vertical map.  But we also know that the left vertical map is injective with
cokernel $\Z_2^{m-1}$, and that the middle vertical map has kernel $\Z_4$ and
cokernel a $\Z_2$-vector space of dimension at most $m$.
Considering all this, we can deduce:
$F_5/F_2\cong \Z_4\oplus \Z_2^{m-1}$, and $E^3_{4,2}=E^{\infty}_{4,2}$.
And $F_5\to \tilde{F_5}$ has cokernel $\Z_2^m$, i.e. the image of $F_5$ in
$\tilde{F_5}$ is exactly $\tilde{F}_2$ which is also the same as all elements
divisible by 2.

Since $\tilde{E}^3_{6,0}$ is free over $\Z$, the short exact sequence for
$\Omega_6^{Spin}(K(\Z^m,2))$ splits. 
Let $\hat{A}(M,x)=\langle e^x \hat{A}(M), [M] \rangle = \langle \frac{4x^3-p_1(M)x}{24}, [M]\rangle$.
Zubr \cite{Zub} proves that there is an injective map 
\begin{gather*}
\Omega_6^{Spin}(K(\Z^m,2))\to \Z^m \oplus \Z^{{m+2}\choose 3}\\
[f:M\to B] \mapsto \left( (\hat{A}(M,f^*(v_i))) , ( \langle f^*(v_iv_jv_k),  [M] \rangle)
\right)
\end{gather*}
with image given by the restrictions 
$$\langle f^*(v_iv_iv_j),  [M] \rangle = \langle f^*(v_iv_jv_j),  [M] \rangle \quad (mod\ 2)$$
for all $i,j$. 
These conditions corresponds exactly to the inclusion of $\tilde{E}^3_{6,0}$ into
$\tilde{E}^2_{6,0}$, and the map to $\Z^m$ is exactly the splitting of the short exact 
sequence to $\tilde{F}_5$.

\begin{thm}\label{TransThm}
The transfer map $\Omega_6^B\to \Omega_6^{\tilde{B}}$ has kernel the torsion subgroup $\Z_4$.
Under Zubr's inclusion map $\Omega_6^{Spin}(K(\Z^m,2))\to \Z^m \oplus \Z^{{m+2}\choose 3}$,
the image of the transfer map corresponds to those classes $[f:M\to B]$ such that
\begin{eqnarray*}
\hat{A}(M,f^*(v_i)) &=& 0 \quad (mod\ 2) \quad \text{for all }i,\\
\langle f^*(v_iv_iv_j), [M] \rangle & = & \langle f^*(v_kv_kv_l), [M] \rangle \ (mod\ 2) 
\quad \text{ for all }i,j,k,l\in\{1,\dots m\}, \\
\langle f^*(v_iv_jv_k), [M] \rangle & = & \frac12 \langle 
f^*(v_i^2v_j + v_iv_j^2+ v_i^2v_k + v_iv_k^2 + v_j^2v_k+ v_jv_k^2), [M] \rangle \ (mod\ 2) \\
& & \quad \text{ for all }i<j<k.
\end{eqnarray*}
\end{thm}
\textbf{Proof:} 
It remains to show that for every $[f:M\to B]\in \Omega_6^B$, we have 
$\hat{A}(\tilde{M},f^*(v_i))= 0 \quad (mod\ 2) \quad \text{for all }i$. 
This would indeed imply that the splitting for the short exact sequence for $\Omega_6^{Spin}(K(\Z^m,2))$ also induces
a splitting for the short exact sequence for $\Omega_6^B$.
By projecting onto the $i$th factor, we may assume $m=1$. 
Recall the fibre bundle $S^2\to \C \Pro \to \HB \Pro$, which restricts to $S^2\to \C \Po^3 \to \HB^1=S^4$.
We use that $\Omega^{Spin}_6(K(\Z,2))\cong \Omega^{Spin}_6(\C \Po^3)\to \Omega^{Spin}_6(S^4)\cong \Z_2$
is given by $$[f:M\to K(\Z,2)] \mapsto \hat{A}(M,f^*v)\ (mod\ 2).$$ 
This follows for example from a comparison of the Atiyah-Hirzebruch spectral sequences.
The map $\C \Pro \times S^\infty \to S^4 \times S^\infty$ is
equivariant with respect to $(\tau,-1)$ on $\C \Pro \times S^\infty$ and $(id,-1)$ on $S^4 \times S^\infty$.
It induces a map between the quotient spaces: $Q_1 \to S^4 \times \R \Pro$.
By naturality of the transfer map it suffices to show that the transfer
$\Omega^{Spin}_6(S^4\times \R \Pro; L)\to \Omega^{Spin}_6(S^4)$ is zero.
But $\Omega^{Spin}_6(S^4\times \R \Pro; L)\cong \Omega_7^{Spin}(S^4\times \R \Pro, S^4)$, 
and via this isomorphism, the transfer map becomes the boundary map in the long exact sequence of the
pair $(S^4\times \R \Pro, S^4)$, see the lemma below. 
But $S^4$ certainly is a retract of $S^4\times \R \Pro$, and so this boundary map is zero.
\qed

\begin{lem}
Given a line bundle $L\to X$, we have the double cover $S(L)\to X$, where $S(L)$ is the corresponding
sphere bundle and $D(L)$ the disk bundle. Under the identification $\Omega_n^{Spin}(X)\cong
\Omega_n^{Spin}(D(L))$ and the Thom isomorphism $\Omega^{Spin}_n(X;L)\to \Omega_{n+1}^{Spin}(D(L),S(L))$,
the transfer $\Omega^{Spin}_n(X;L)\to \Omega^{Spin}_n(S(L))$ corresponds to the boundary map
$\Omega_{n+1}^{Spin}(D(L),S(L))\to \Omega_n^{Spin}(S(L))$ in the long exact sequence of the pair $(D(L),S(L))$.
\end{lem}
\textbf{Proof:}
Given $[f:(M^{n+1},\partial M)\to (D(L),S(L))]\in \Omega_{n+1}^{Spin}(D(L),S(L))$, we may assume
$f$ is transverse to the zero section $X$ of the line bundle. Then the corresponding element in
$\Omega^{Spin}_n(X;L)$ is given by $[f^{-1}(X)\to X]$, and its image under the transfer map is
represented by $f^{-1}(X)\times_{X}S(L)\to S(L)$. But this is spin bordant to $\partial M \to S(L)$:
Let $p:[0,1] \times S(L)\to D(L)$ be given by $(t,v)\mapsto tv$. 
Then $M\times_{D(L)} ([0,1]\times S(L))\to S(L)$ defines the spin bordism we need.
The geometric interpretation is that we cut $M$ along $f^{-1}(X)$.
\qed

\subsection{Identifying the result of the construction}

We want to know which are the conjugation manifolds $X$ we obtain by our construction. 
So we should (by Wall's result) compute the first Pontrjagin class of $X$ 
and the triple cup products of degree 2 cohomology classes.

By Zubr's result it follows that if we equip our conjugation space $X$ with a map to 
$K(\Z^{m-1},2)$ which is an isomorphism on $\pi_2$, then it suffices to know the spin 
bordism class of $X$.

The process of starting with a manifold $W$ with normal $B$-structure and boundary $M\times
\R \Po^2$ and taking the union of the double cover $V$ with $M\times D^3$ along the common boundary
$M\times S^2$ as described earlier makes sense for any such $W$, not only those
$W$ which produce
conjugation spaces. It can be interpreted as a map of sets of bordism classes
$$\Omega_6^B(\partial=M\times \R \Po^2) \to \Omega_6^{Spin}(K(\Z^m,2), \partial = M\times S^2)
 \to \Omega_6^{Spin}(K(\Z^{m-1},2)).$$
We still have to make precise how to obtain the map $X\to K(\Z^{m-1},2)$ from $V\to K(\Z^m,2)$:
Let $p:K(\Z^m,2)\to K(\Z^{m-1},2)$ be the map induced by the group homomorphism
$\Z^m\to \Z^{m-1}$ defined by $e_i\mapsto e_i$ for $i<m$ and $e_m\mapsto - \sum e_i$.
The use of this map $p$ corresponds to the choice of the basis
given by $e_1,\dots e_{m-1}$ for 
$$\langle e_i \rangle / \langle \sum e_i \rangle \cong H_2(K(\Z^m,2))/Im H_2(S^2)$$
which is isomorphic to $H_2(V)/Ker(H_2(V)\to H_2(X))\cong H_2(X)$
in the case where we produce a conjugation space $X$.
We see that in integral cohomology $p^*$ sends $v_i\mapsto v_i-v_m$.
We obtain a map $V\to K(\Z^{m-1},2)$ as the composition of $V\to K(\Z^m,2)$ with $p$.
The restriction of $V\to K(\Z^{m-1},2)$ to the boundary $M\times S^2$ is nullhomotopic,
which implies that the corresponding element in $H^2(V;\Z^{m-1})$ comes from 
$H^2(V,M\times S^2;\Z^{m-1})\cong H^2(X,M;\Z^{m-1})$. Since $H^2(X)\to H^2(V)$ is injective, this means that the
element comes from a uniquely determined $H^2(X,\Z^{m-1})$, or in other words that
the map $V\to K(\Z^{m-1},2)$ extends uniquely up to homotopy to $X\to K(\Z^{m-1},2)$.

In the above we use bordism classes of manifolds with a fixed boundary.
These bordism sets of manifolds with boundary 
are torsors $\Omega_6^B(\partial=M\times \R \Po^2)$ for the bordism group $\Omega_6^B$, and
$\Omega_6^{\tilde{B}}(\partial=M\times S^2)$ for the bordism group
$\Omega_6^{\tilde{B}}=\Omega_6^{Spin}(K(\Z^m,2))$.

Moreover, the map   
$$
\Omega_6^B(\partial=M\times \R \Po^2) \to \Omega_6^{Spin}(K(\Z^{m-1},2))
$$
is equivariant with respect to the group homomorphism given as the composition 
$$
\Omega_6^B \stackrel{tr}{\to} \Omega_6^{Spin}(K(\Z^m,2))
\stackrel{p_*}\to \Omega_6^{Spin}(K(\Z^{m-1},2)).
$$

To compute the group homomorphism, we use the injections 
$\Omega_6^{Spin}(K(\Z^{r},2))\to \Z^r \oplus H_6(K(\Z^r,2))$
given by the $\hat{A}$-genera and triple cup products.
Let us denote the generators of $\Z^r$ by $d_i$ and the generators of $H_6(K(\Z^r,2))$ by
$e_{ijk}$ as before. 

The image of $[M\stackrel{f}\to K(\Z^m,2)]$ under $p_*$ is $[M\stackrel{p\circ f}\to K(\Z^{m-1},2)]$. 
In order to write it as a linear combination $\sum_i \lambda_i d_i + \sum_{i\le j \le k} \mu_{ijk} e_{ijk}$ we have to 
compute the coefficients 
\begin{eqnarray*}
\lambda_i & = & \hat{A}(M,(p\circ f)^*(v_i))\\
 & = & \hat{A}(M,f^*(v_i-v_m))\\
 & = & \hat{A}(M,f^*v_i)- \hat{A}(M,f^*v_m) + \frac12 \langle f^*(-v_i^2v_m +v_iv_m^2) , [M] \rangle 
 \qquad \text{and}\\ 
\mu_{ijk} & = & \langle (p\circ f)^*(v_iv_jv_k),[M]\rangle \\
 & = & \langle f^*((v_i-v_m)(v_j-v_m)(v_k-v_m)),[M]\rangle \\
 & = & \langle f^*(v_iv_jv_k - v_iv_jv_m - v_iv_kv_m - v_jv_kv_m + v_iv_m^2+v_jv_m^2+v_kv_m^2-v_m^3),[M]\rangle .
\end{eqnarray*}
  
Thus we obtain: 
\begin{eqnarray*}
d_i & \mapsto & d_i \quad \text{for }i < m,\\
d_m & \mapsto & -\sum_i d_i, \\
e_{ijk} & \mapsto & e_{ijk} \quad \text{for }i\le j\le k < m,\\
e_{ijm} & \mapsto & -\sum_{k\ne i,j} e_{ijk} - 2 e_{iij} - 2 e_{ijj} \quad \text{for }i<j<m,\\
2 e_{iim} & \mapsto & - d_i - 2 \sum_{j\ne i} e_{iij} - 6 e_{iii} \quad \text{for }i< m,\\
e_{iim} + e_{imm} & \mapsto & \sum_{\stackrel{j\le k}{j,k\ne i}} e_{ijk} +  \sum_{j\ne i} e_{iij} 
\quad \text{for }i< m,\\
e_{mmm} & \mapsto & -\sum_{i\le j\le k} e_{ijk}.  
\end{eqnarray*}

This implies that under $p_*:\Omega_6^{Spin}(K(\Z^m,2))\to \Omega_6^{Spin}(K(\Z^{m-1},2))$,
the generators of the image of the transfer map are mapped in the following way:
\begin{eqnarray*}
2 d_i & \mapsto & 2 d_i \quad \text{for }i < m,\\
2 d_m & \mapsto & -\sum_i 2 d_i,\\
\sum_{i\le j\le k} e_{ijk} & \mapsto & 0,\\ 
2e_{iii} & \mapsto & 2 e_{iii} \quad \text{ for }i<m, \\ 
\sum_{j=1}^m 2e_{iij} & \mapsto & -d_i - 6e_{iii} \quad \text{ for }i<m, \\
2e_{imm}+\sum_{k\ne i, k\ne m} e_{imk} & \mapsto & d_i + 6e_{iii} +2 \sum_{j\ne i} e_{iij} \quad \text{ for }i<m, \\
2e_{ijj}+\sum_{k\ne i, k\ne j} e_{ijk} & \mapsto & -2 e_{iij} \quad \text{ for }i\ne j\text{ and
}i,j<m, \\
2e_{ijk} & \mapsto & 2e_{ijk} \text{ for }i<j<k.
\end{eqnarray*}

Thus the image of $\Omega_6^B \to \Omega_6^{Spin}(K(\Z^{m-1},2))$ has generators
$d_i\text{ for }1\le i \le m-1$, and $2 e_{ijk}\text{ for } 1\le i\le j\le k \le m-1$.

\bigskip

We know that the construction produces manifolds $X$ such that the $\Z_2$-cohomology ring 
of $X$ is isomorphic to the $\Z_2$-cohomology ring of $M$ with all degrees multiplied by 2.
But our calculation shows that by choosing different $W$ in different bordism classes
producing different conjugation 6-manifolds $X$ with the same fixed point set $M$,
we can modify the bordism class of $X$ by $d_i$ and $2e_{ijk}$, which means that we can produce
conjugations with fixed point set $M$ on any simply-connected spin 6-manifold with free
cohomology, which has $\Z_2$ cohomology ring isomorphic to the $\Z_2$ cohomology ring of $M$
using an isomorphism dividing all degrees by 2. This proves theorem \ref{nthm}.

Finally, we prove corollary \ref{mcor}.
Postnikov  proved that those trilinear forms $t$ on 
finite-dimensional $\Z_2$-vector spaces $V$ with values in $\Z_2$ which occur as the triple cup
products of one-dimensional $\Z_2$-cohomology classes of a closed oriented 3-manifold
are exactly those which satisfy $t(v,v,w)=t(v,w,w)$ for all $v,w\in V$ \cite{Pos}.
If $H_1(M;\Z)$ contains no summand $\Z_2$, then we have the additional condition $t(v,v,w)=0$ for all $v,w\in V$
(since $H^2(M;\Z)$ contains no summand $\Z_2$, and so $v^2=Sq^1(v)=0$ for all $v\in H^1(M;\Z_2)$). 
Sullivan \cite{Sul} constructed closed oriented 3-manifolds with free integral cohomology
realizing all trilinear forms $t$ satisfying this additional condition.

(Matthias Kreck explained to me how to prove the realizability in a different way: 
The additional condition means exactly that the trilinear form comes from 
a $\Z$-valued trilinear form on a free $\Z$-module. 
Trilinear forms on $\Z^m$ are in bijection with $H_3(K(\Z^m,1))$. 
Now $\tilde{\Omega}_3^{Spin}(K(\Z^m,1))\to H_3(K(\Z^m,1))$ 
is surjective, as one can see from the Atiyah-Hirzebruch spectral sequence.
By surgery below the middle dimension, every bordism class contains 2-connected 
maps $M\to K(\Z^m,1)$, so that $H^1(M;\Z)\cong \Z^m$, with the desired trilinear 
form given by triple cup products.)

{\textsc{Max-Planck-Institut f\"ur Mathematik, Bonn, Vivatsgasse 7, 53111 Bonn, Germany}}\\
E-mail: \nolinkurl{olber@mpim-bonn.mpg.de}

\end{document}